\documentclass{elsarticle}

\usepackage{amssymb}
\usepackage{amsmath}
\newtheorem{theorem}{Theorem}
\newtheorem{proposition}{Proposition}[section]

\begin{document} 

\begin{frontmatter} 
  
\title{On existence of a Morse energy function for topological flows with finite chain recurrent sets}


\author[mtvadd]{Timur~V.~Medvedev\corref{mycorrespondingauthor}}
\cortext[mycorrespondingauthor]{Corresponding author} 
\address[mtvadd]{Laboratory of Algorithms and Technologies for Network Analysis; National Research University Higher School of Economics; 136 Rodionova Street, Niznhy Novgorod, Russia}
\ead{mtv2001@mail.ru}

\author[opadd]{Olga~V.~Pochinka}
\address[opadd]{National Research University Higher School of Economics; 25/12 Bolshaya Pecherckaya Street, Niznhy Novgorod, Russia}
\ead{olga-pochinka@yandex.ru}

\author[szadd]{Svetlana~Kh.~Zinina}
\address[szadd]{Ogarev Mordovia State University;111v-2, Bolshevistskaya street, Saransk, Russia}
\ead{suddenbee@gmail.com}

\begin{abstract}
We prove the existence of a continuous Morse energy function for an arbitrary topological flow with finite hyperbolic (in topological sense) chain recurrent set on a topological manifold of any dimension. This result is a partial solution of the Morse problem of existence of continuous Morse functions on any topological manifolds. Namely, we prove that a topological manifold admits a continuous Morse function if it admits a topological flow with finite hyperbolic chain recurrent set.
\end{abstract}

\end{frontmatter}

\section{Introduction and main results}
\label{sec1}

It is well known that for dimensions from 4 on there are topological manifolds admitting no smooth structure, therefore, dynamical systems as well as functions on such manifolds may only be considered as topological and continuous, respectively. Nevertheless, these systems and functions have the same properties as the smooth ones and they are closely related to the topology of the ambient manifold. In this paper we consider a well known {\it Lyapunov function} for dynamical systems, that is a continuous function which is constant on each chain component and decreases along orbits outside the chain recurrent set. More precisely, we consider an {\it energy function} whose set of critical points coincides with the chain recurrent set of the system.

Definitions.

Let $M^n$ be a closed $n$-manifold with a metric. {\it A topological flow} on $M^n$ is a family of homeomorphisms $f^t:M^n\to M^n$ each of which is continuous on $t\in\mathbb R$ and satisfies 
\begin{enumerate}
	\item $f^0(x)=x$ for any $x\in M^n$;
	\item $f^t(f^s(x))=f^{t+s}(x)$ for any $s,t\in \mathbb{R}$, $x\in M^n$. 
\end{enumerate}

The {\it trajectory} or the {\it orbit} of a point $x\in M^n$ with respect to the flow $f^t$ is the set $\mathcal O_x=\{f^t(x), t\in \mathbb{R}\}$. The trajectories are oriented with respect to the parameter $t$. Any two trajectories either coincide or they do not intersect one another, therefore, the phase space is the union of pairwise disjoint trajectories. There are three types of trajectories:

\begin{enumerate}
	\item a {\it fixed point} $\mathcal O_x=\{x\}$.
	\item a {\it periodic trajectory} $\mathcal O_x$ for which there exists such $per(x)>0$ that $f^{per(x)}(x)=x$ but $f^t(x)\neq x$ for every $0<t<per(x)$. The number $per(x)$ is called the {\it period} of $\mathcal O_x$ and it is independent of the choice of $x\in\mathcal O_x $.
	\item a {\it regular} trajectory $\mathcal O_x$ is a trajectory that is neither a fixed point nor a periodic trajectory. Every regular trajectory is homeomorphic to the line.
\end{enumerate}

An {\it $\varepsilon$-chain of length $T$} connecting a point $x$ to a point $y$ with respect to the flow $f^t$ is a sequence of points $x=x_0,\dots,x_n=y$ for which there is a sequence $t_1,\dots,t_n$ such that $d(f^{t_i}(x_{i-1}),x_{i})<\varepsilon$,  $t_i\geq 1$ for $1\leq i
\leq n$ and $t_1+\dots+t_n=T$.

A point $x\in M^n$ is {\it chain recurrent} for the flow $f^t$ if for every $\varepsilon>0$ there is $T$ and there is an $\varepsilon$-chain of length $T$ connecting $x$ to itself. The set of all chain recurrent points of $f^t$ is the {\it chain recurrent set} denoted by $\mathcal R_{f^t}$, its connected components being the  {\it chain components}. The set $\mathcal R_{f^t}$ is $f^t$-invariant, i.e. it is composed of the $f^t$ trajectories which are called {\it chain recurrent}. Fixed points and periodic orbits are chain recurrent. 

A {\it Lyapunov function} for a flow is a continuous function which decreases along the orbits outside the chain recurrent set and which is constant on each chain component. It follows from the results of C.~Conley \cite{Con1978} that such a function exists for every flow defined by a continuous vector field (this fact is known as ``The Fundamental Theorem of Dynamical Systems''). From the results of W.~Wilson \cite{wilson1969smoothing} it follows that for every such a flow there is a smooth Lyapunov function whose set of critical points coincide with the chain recurrent set. Such Lyapunov function is called the {\it energy function}. From the results of S.~Smale \cite{S12} and K.~Meyer \cite{Me1968} it follows that any gradient-like flow on a manifold has a Morse energy function. 
	
In this paper we develop the ideas of \cite{S12} and \cite{Me1968}. In section \ref{sec2} we introduce the notions of a topologically hyperbolic fixed point and a continuous Morse function on a manifold $M^n$. For the class $G(M^n)$ of topological flows with finite (therefore, composed of fixed points) hyperbolic chain recurrent set on $M^n$ we prove
	
	\begin{theorem}\label{mainr} Every flow $f^t\in G(M^n)$ admits a Morse energy function.
	\end{theorem}
	
	\section{Dynamics of flows of $G(M^n)$}
	\label{sec2} 
	Let $f^t$ be a topological flow on a closed manifold $M^n$. 
	
	Two flows $f^t:M^n\to M^n$, $g^t:M^n\to M^n$ are said to be {\it topologically equivalent} if there is a homeomorphism $h:M^n\to M^n$ such that $g^th=hf^t$ for every $t\in\mathbb R$. $h$ is called the {\it conjugating} homeomorphism.
	
	Let the model flow in a neighborhood of a fixed point be the linear flow $a^t_\lambda:\mathbb R^n\to\mathbb R^n, \lambda \in \{0,1, \ldots , n \}$ defined by
		  $$ a ^ t_ {\lambda} (x_1, \ldots , x_ {\lambda}, x_ {\lambda + 1}, \ldots , x_n) = (2 ^ tx_1, \ldots , 2 ^ tx_ {\lambda }, 2 ^ {- t} x_ {\lambda + 1}, \ldots , 2 ^ {- t} x_n). $$ Let $$ {{E} ^ s_\lambda} = \{(x_1, \ldots , x_n) \in \mathbb {R} ^ n: x_1 = \dots = x_ {\lambda} = 0 \},$$ $$ {{E} ^ u_\lambda} = \{(x_1, \ldots , x_n) \in \mathbb {R} ^ n: x_ {\lambda + 1} = \dots = x_ {n} = 0 \}. $$
	
	A fixed point $p$ is {\it topologically hyperbolic} if there are a neighborhood $U_p\subset M^n$ of $p$, a number $\lambda \in \{0,1, \ldots , n \}$ and a homeomorphism $h_p:U_p\to \mathbb{R}^n$ such that $h_pf^t|_{U_p}=a^t_{\lambda_p}h_p|_{U_p}$ whenever both the right and the left side are defined.
The {\it locally invariant manifolds} of a fixed point $p$ are the sets $h_p^{-1}(E^s_{\lambda_p}), h_p^{-1}(E^u_{\lambda_p})$. The sets $$W^s_p=\bigcup \limits_{t\in \mathbb{R}}f^t(h_p^{-1}(E^s_{\lambda_p})),\,W^u_p=\bigcup \limits_{t\in \mathbb{R}}f^t(h_p^{-1}(E^u_{\lambda_p}))$$ are called the {\it stable} and the {\it unstable} invariant manifolds of $p$, respectively. If follows from the definition that $W^s_p=\{y\in M^n: \lim\limits_{t\to +\infty} f^{t}(y)= p\}, W^u_p=\{y\in M^n: \lim\limits_{t\to +\infty} f^{-t}(y)= p\}$ and $W^u_p\cap W^u_q=\emptyset$ ($W^s_{p}\cap W^s_{q}=\emptyset$) for any two hyperbolic fixed points $p,q$. Moreover, there is an injective immersion $J:\mathbb{R}^{\lambda_p}\to M^n$ such that $W^u_p=J(\mathbb{R}^{{\lambda_p}})$\footnote{A map $J:\mathbb{R}^{m}\to M^n$ is an {\it immersion} if for every point $x\in \mathbb{R}^{m}$ there is such a neighborhood $U_x\subset \mathbb{R}^{m}$ that the restriction $J|_{U_x}$ of the map $J$ to $U_x$ is a homeomorphism.}.

	We say the number $\lambda_p$ to be the {\it index} of the fixed hyperbolic point $p$. We say points of indexes $n$ and $0$ to be {\it sources} and {\it sinks}, respectively. A point $p$ with  $\lambda_p\in\{1,\cdots,n-1\}$ is said to be a {\it saddle}. 
	
Recall that the flows of $G(M^n)$ are topological flows on $M^n$ with finite hyperbolic chain recurrent set, therefore, the points of the chain recurrent set are the fixed points. The dynamics of these flows are similar to that of the gradient-like flows in the following sense.  Analogously to the Smale's order we introduce the relation on the set of fixed points of $f^t\in G(M^n)$ by $$p\prec q\iff W^s_p\cap W^u_q\neq\emptyset.$$ Since the chain recurrent set of $f^t$ is finite this relation can be extended to a total order relation on $\mathcal R_{f^t}$. From now on let the fixed points of $f^t$ be enumerated according to this order: $$p_1\prec\dots\prec p_k.$$ Assume without loss of generality that in this order any sink precedes any saddle and any saddle precedes any source. 
	
Using the methods analogous to that of \cite{GrGuMePo2017} one can prove the following theorem which describes the embedding and the asymptotic behavior of the invariant manifolds of the fixed points.
	
	\begin{theorem}\label{dynam} Let $f^t\in G(M^n)$. Then
		\begin{enumerate}
			\item  $M^n=\bigcup\limits_{i=1}^kW^u_{p_i}=\bigcup\limits_{i=1}^kW^s_{p_i}$;
			\item  $W^u_{p_i}\,\,(W^s_{p_i})$ is a topological submanifold of $M^n$ and it is homeomorphic to $\mathbb R^{\lambda_{p_i}}\,(\mathbb R^{n-\lambda_{p_i}})$;
			\item
			$cl(W^u_{p_i})\setminus
			(W^u_{p_i}\cup p_i)\subset\bigcup\limits_{j=1}^{i-1}W^u_{p_j}\,\,\, (cl(W^s_{p_i})\setminus
			(W^s_{p_i}\cup p_i)\subset\bigcup\limits_{j=i+1}^{k}W^s_{p_j})$.
		\end{enumerate}
	\end{theorem}

\section{Continuous Morse function}
	\label{sec3}
	Following M.~Morse \cite{morse1959topologically} we now introduce a continuous Morse function on  $M^n$.
	
	Let $\varphi:M^n\to\mathbb R$ be a continuous function with real values. A point $p\in M^n$ is said to be {\it regular} if there is a neighborhood $V_p\subset M^n$ of $p$ and there is a homeomorphism $\phi_p:y\in V_p\mapsto \phi_p(y)=(x_1(y),\cdots, x_n(y))\in\mathbb R^n$ such that 
	$$x_i(p) = 0, i\in\{1,\cdots, n\},\,\,\varphi(y)=\varphi(p)+x_n(y),\,y\in V_p.$$ Otherwise the point $p$ is called {\it critical}. Denote by $Cr_\varphi$ the set of critical points of $\varphi$. If the coordinates $x_i,i\in\{1,\cdots, n\}$ of a critical point $p$ are such that there is $\nu_p\in\{0,\cdots,n\}$ such that  $$\varphi(y)=\varphi(p)-\sum\limits_{i=1}^{\nu_p}x^2_i(y)+\sum\limits_{i=\nu_p+1}^nx^2_i(y),\,\,y\in V_p,$$ then $p$ is called a {\it non-degenerate critical point of index $\nu_p$}. A Morse function whose every critical point is non-generate is called a {\it continuous Morse function}.  
	
	From \cite{morse1959topologically} it follows that the classical Morse inequities for the number of critical points of index $\nu$ and the $\nu$-th Betti number are true for a continuous Morse function.   
	
	Proposition \ref{action_on_W} stating the existence of a local Morse energy function follows from hyperbolicity of the fixed point $p$.
	
	\begin{proposition}\label{action_on_W} Let $p$ be a fixed point of index $\lambda_p$ of a flow $f^t\in G(M^n)$ and let  $h_p:y\in U_p\mapsto h_p(y)=(x_1(y),\cdots, x_n(y))\in\mathbb{R}^n$ be a homeomorphism conjugating $f^t$ in a neighborhood ${U_p}$ of $p$ to the linear flow $a^t_{\lambda_p}$. Then for every number $c\in\mathbb R$ the function $$\varphi_{p,c}(y)=c-\sum\limits_{i=1}^{\lambda_p}x^2_i(y)+\sum\limits_{i=\lambda_p+1}^nx_i^2(y),\,\,y\in U_p$$ is the local Morse energy function for $f^t$ in the neighborhood of $p$.
	\end{proposition}

\section{Construction of an energy function for $G(M^n)$ flows}\label{sec4}

In this section we prove Theorem \ref{mainr}, that is for every topological flow $f^t$ of $G(M^n)$ we construct a continuous Morse function $\varphi:M^n\to\mathbb R$ with the properties
\begin{enumerate}	
\item $\varphi(f^t(x))<\varphi(x)$ for every $x\in(M^n\setminus\mathcal R_{f^t})$ and every $t>0$;
	
\item $Cr_\varphi=\mathcal R_{f^t}$ and $\nu_{p_i}=\lambda_{p_i}$ for any fixed point $p_i$.
\end{enumerate}
	
For each $i\in\{1,\dots,k-1\}$ let $$A_i=\bigcup\limits_{j=1}^iW^u_{p_j},\,\,R_i=\bigcup\limits_{j=i+1}^iW^s_{p_j}.$$ By induction on $i\in\{1,\dots,k-1\}$ we are going to construct a neighborhood $U_i$ of the set $A_i$ and a Morse function $\varphi_i:U_i\to[1,i+1/3]$ with the following properties:
\begin{itemize}
\item $U_{i}=\varphi_i^{-1}([1,i+1/3])$ is a closed compact topological $n$-submanifold with the boundary $\partial U_i=\varphi_i^{-1}(i+1/3)$ such that $A_i\subset\mathop\mathrm{int}\,U_i,\,U_i\cap R_i=\emptyset$.

\item $\varphi_i$ is an energy function for the flow $f^t|_{U_i}$ which coincides with $\varphi_{p_i,i}$ in some neighborhood of $p_i$. 
\end{itemize}

The last step of the construction that is the construction of a neighborhood $U_{k}$ and a function $\varphi_k$ when $\varphi_{k-1}$ is already constructed on the neighborhood $U_{k-1}$ we consider separately. According to Theorem \ref{dynam} the neighborhood $U_k$ is the entire manifold $M^n$ and the function $\varphi_k$ is the desired function $\varphi$.

Denote by\\
$\mathbb{D}^n=\{(x_1,\dots,x_{n})\in\mathbb{R}^{n}~:~
\sum\limits_{i=1}^{n}x_i^2\leq 1\}$ the standard $n$-disk ($n$-ball), $\mathbb D^0=\{0\}$ and by\\
$\mathbb{S}^{n-1}=\{(x_1,\dots,x_{n})\in\mathbb{R}^{n}~:~\sum
\limits_{i=1}^{n}x_i^2=1\}$ the standard $(n-1)$-sphere, $\mathbb S^{0}=\{-1,1\}$.
	
{\it Construction for $i=1$.} By Theorem \ref{dynam} the point $p_1$ is a sink. By Proposition \ref{action_on_W} there is a local energy function $\varphi_{p_1,1}$ in some neighborhood of $p_1$. Let $U_1=\varphi^{-1}_{p_1,1}([1,4/3])$ and $\varphi_1=\varphi_{p_1,1}|_{U_1}$ and this concludes the proof for this step.
	
{\it Induction step.} Suppose the desired $\varphi_{i-1}:U_{i-1}\to[1,i-2/3]$ is constructed. Now we are going to construct $\varphi_{i}:U_{i}\to[1,i+1/3]$. Consider three cases:  $p_i$ is a) a sink; b) a saddle; c) a source.
	
\begin{figure}[ht]
	\centering
	\includegraphics[width=0.7\linewidth]{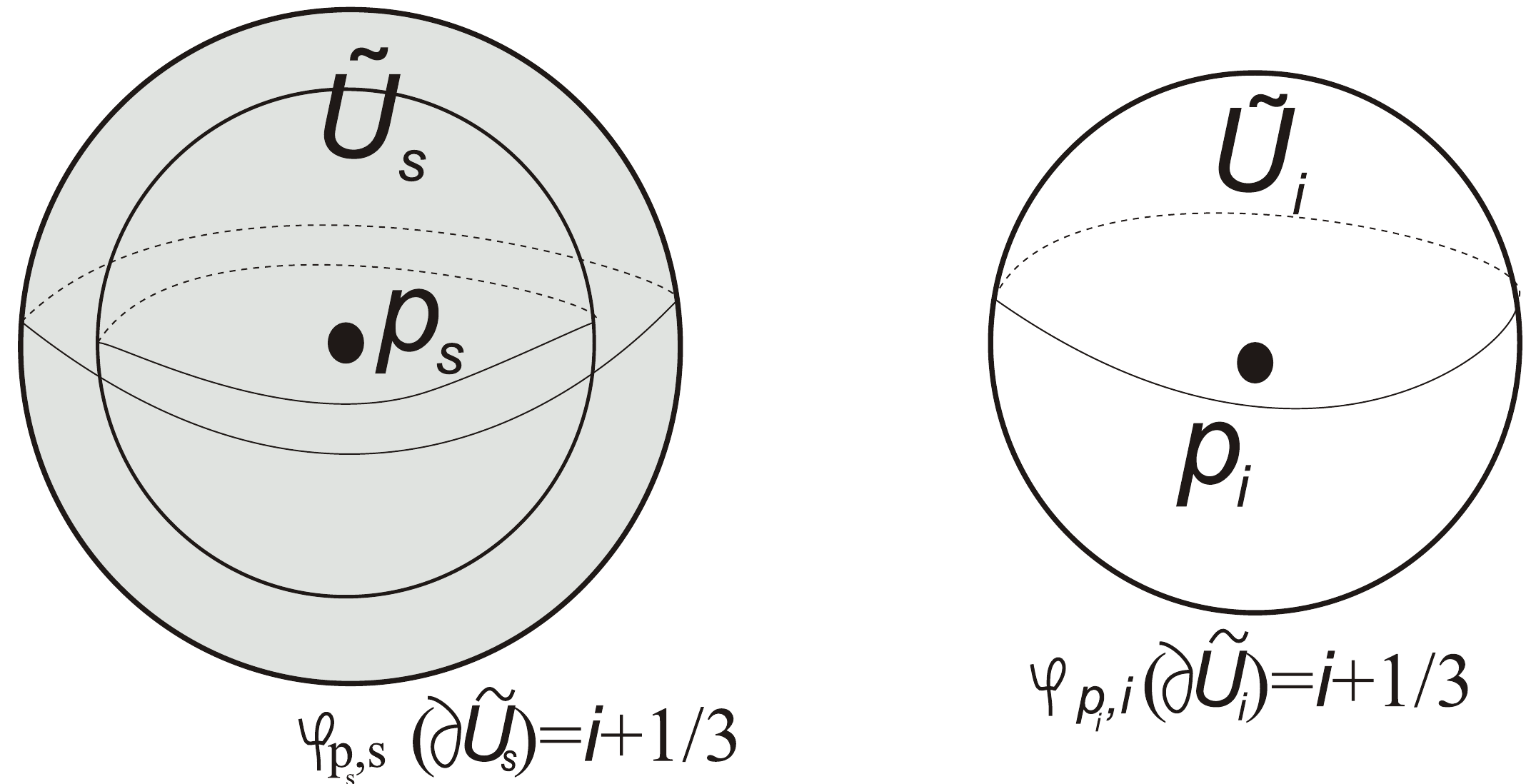}
	\caption{Induction step if $p_i$ is a sink}
	\label{fig2}
\end{figure}

{\it a) The point $p_{i}$ is a sink.} Analogous to the step $i=1$ from Proposition \ref{action_on_W} it follows that  for each $s\in\{1,\dots,i\}$ there is a local energy function ${\varphi}_{p_{s},s}$ in some neighborhood of $p_{s}$. Let $\tilde{U}_{s}={\varphi}^{-1}_{p_{s},s}([s,i+1/3])$, $U_i=\bigcup\limits_{s=1}^{i}\tilde U_s$ and define the desired function  $\varphi_i$ by
	$$\varphi_{i}(x)=\varphi_{p_{s},s}(x),  x\in \tilde U_{s},$$ (see Figure \ref{fig2}).
	
\begin{figure}[ht]
\centering
\includegraphics[width=1.0\linewidth]{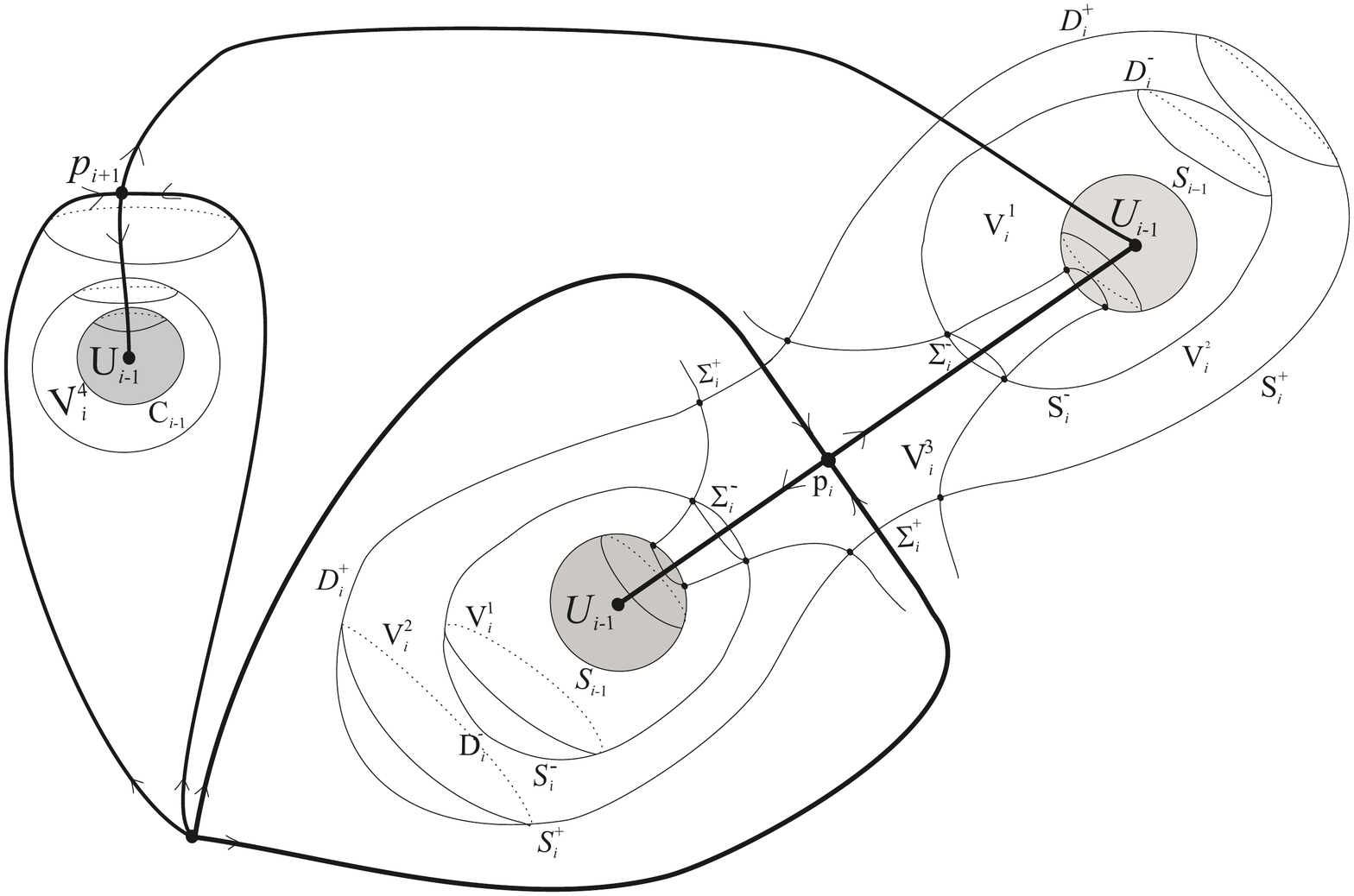}
\caption{Induction step if $p_i$ is a saddle}
\label{fig3}
\end{figure}

{\it b) The point $p_{i}$ is a saddle with Morse index $\lambda_{p_i}$.} From Proposition \ref{action_on_W} it follows that there is a local energy function ${\varphi}_{p_{i},i}$ in some neighborhood of $p_{i}$. Let  $L=\{(x_1,\dots,x_n)\in\mathbb R^n:\,x_{\lambda_{p_i}+1}^2+\cdots+x^2_n\leq 1/4\}$ and  $\Sigma^\pm=\{(x_1,\dots,x_n)\in\mathbb R^n:\,-x_1^2-\cdots-x_{\lambda_{p_i}}^2+x_{\lambda_{p_i}+1}^2+\cdots+x^2_n=\pm 1/3\}$. By construction the set $L\subset\mathbb R^n$ is homeomorphic to $\mathbb R^{\lambda_{p_i}}\times \mathbb D^{n-\lambda_{p_i} }$ and the set  $\Sigma^-\,(\Sigma^+)$ is homeomorphic to $\mathbb{S}^{\lambda_{p_i}-1}\times \mathbb R^{n-\lambda_{p_i}}\,(\mathbb{S}^{n-\lambda_{p_i}-1}\times \mathbb R^{\lambda_{p_i}})$. Let   $$\Sigma_i^\pm={\varphi}^{-1}_{p_{i},i}(i\pm1/3),\,\,Q_i=\Sigma^-_i\cap W^u_{p_i},\,\,d_i=\Sigma^-_i\cap{\varphi}^{-1}_{p_{i},i}(L).$$ Then these sets satisfy: 
\begin{itemize}
\item the set $\Sigma_i^\pm\subset U_{p_i}$ is the image of $\Sigma^\pm$ by the homeomorphism $h_{p_i}^{-1}$, therefore,  $\Sigma_i^-\,(\Sigma_i^+)$ is homeomorphic to $\mathbb{S}^{\lambda_{p_i}-1}\times \mathbb R^{n-\lambda_{p_i}}\,(\mathbb{S}^{n-\lambda_{p_i}-1}\times \mathbb R^{\lambda_{p_i}})$;

\item the set $Q_i\subset \Sigma^-_i$ is the image of $\mathbb R^{\lambda_{p_i}}\cap \Sigma^-$ by the homeomorphism $h_{p_i}^{-1}$ and, therefore, it is homeomorphic to $\mathbb{S}^{\lambda_{p_i}-1}$;

\item the set $d_i\subset\Sigma^-$ is the image of $L\cap \Sigma^-$ by the homeomorphism $h_{p_i}^{-1}$ and, therefore, it is homeomorphic to $\mathbb{S}^{\lambda_{p_i}-1}\times\mathbb D^{n-\lambda_{p_i} }$ (see Figure \ref{fig3}).  
\end{itemize}

Without loss of generality assume $d_i\cap U_{i-1}=\emptyset$ (to satisfy this condition one changes values of the function ${\varphi}_{p_{i},i}$). Let $c_{i}=\partial d_{i}$ and $R_{i}=\bigcup\limits_{x\in c_{i}}\mathcal O_x$. Then  $c_{i}$ is homeomorphic to $\mathbb{S}^{\lambda_{p_i}-1}\times \mathbb{S}^{n-\lambda_{p_i}-1}$ and $R_{i}$ is homeomorphic to $\mathbb{S}^{\lambda_{p_i}-1}\times \mathbb{S}^{n-\lambda_{p_i}-1}\times\mathbb{R}.$ According to the order relation every trajectory $\mathcal O_x,x\in d_i$ intersects $\partial U_{i-1}$ at a single point. Let $$d_{i-1}=\bigcup\limits_{x\in d_i}(\mathcal O_x\cap\partial U_{i-1}),\,c_{i-1}=\bigcup\limits_{x\in c_i}(\mathcal O_x\cap\partial U_{i-1}).$$ Then $f^{-t_1(x)}(x)=\mathcal O_x\cap d_i,x\in d_{i-1}$ defines a positive continuous function $t_1:d_{i-1}\to\mathbb R$ mapping $x\in d_{i-1}$ to $t_1(x)$.
	
By the induction hypothesis each connected component of $\partial U_{i-1}$ is a closed topological $(n-1)$-manifold. Denote by  $S_{i-1}$ the union of the connected components of $\partial U_{i-1}$ intersecting the set $d_{i-1}$. Let $B_{i-1}=S_{i-1}\setminus d_{i-1}$. By construction $B_{i-1}$ is a compact $(n-1)$-submanifold with the boundary $c_{i-1}$ and it is homeomorphic to $\mathbb{S}^{\lambda_{p_i}-1}\times \mathbb{S}^{n-\lambda_{p_i}-1}$. From the {\it collaring theorem} (see, for example, Theorem 2.1, p.152 \cite{hir94}) there exists an embedding $\xi: c_{i-1}\times [0,1]\to B_{i-1}$ for which  $\xi(c_{i-1}\times\{0\})=c_{i-1}.$ Let $E_{i-1}=\xi(c_{i-1}\times [0,1])$.	

Then for every point $y\in E_{i-1}$ there exists a unique pair of points $x_y\in c_{i-1}$ and $t_y\in[0,1]$ such that $y=\xi (x_y,t_y)$. On the other hand since $x_y\in d_{i-1}$ there exists a corresponding $t_1(x_y)$ for every $x_y$. 
Extend the function $t_1$ to the positive function $L_1:E_{i-1}\to\mathbb R$ by $L_{1}(y)=t_1(x_y)+t_y(1-t_1(x_y)).$ Define the function $T_1:S_{i-1}\to\mathbb R$ by
$$T_{1}(x)=\begin{cases}
{t}_{1}(x),  x\in d_{i-1}, \cr
L_{1}(x),  x\in E_{i-1}, \cr
1, x\in (B_{i-1}\setminus E_{i-1}).
\end{cases}$$
Let $$S^-_{i}=\bigcup\limits_{x\in S_{i-1}}f^{-T_1(x)}(x),\,\,\,\,V^1_{i}=\bigcup\limits_{x\in S_{i-1}}\left(\bigcup\limits_{t\in[0,T_1(x)]}f^{-t}(x)\right).$$
Thus every point $y\in V_i^1$ is uniquely represented as $y=f^{-t_y}(x_y)$ where $x_y\in S_{i-1}$ and $t_y\in[0,T_1(x_y)]$. Define the function $\varphi_{_{V_i^1}}:V_i^1\to[i-2/3,i-1/3]$ by $$\varphi_{_{V_i^1}}(y)=i-\frac{2-t_y/T_1(x_y)}{3}.$$ 
	
Let $D^-_{i}=S^-_{i}\setminus d_{i}$. By construction $D^-_i$ is the compact $(n-1)$-submanifold with the boundary $c_{i}$. Denote by $t_2>0$ the time moment such that $f^{-t_2}(x)\in\Sigma^+_i$ for $x\in c_{i}$. For every $t\in[0,t_2]$ let $\psi(t)={\varphi}_{p_{i},i}(f^{-t}(c_i))$. Let $D^+_{i}=\bigcup\limits_{x\in D^-_{i}}f^{-t_2}(x)$ and denote by $V_i^2\subset M^n$ the compact set bounded by the compact $(n-1)$-submanifolds $D^-_{i}$, $D^+_{i}$, $R_{i}$. Every point $y\in V_i^2$ is uniquely represented as $y=f^{-t_y}(x_y)$ where $x_y\in D^-_{i}$ and  $t_y\in[0,t_2]$. Define the function $\varphi_{_{V_i^2}}:V_i^2\to[i-1/3,i+1/3]$ by $$\tilde\varphi_{_{V_i^2}}(y)=\psi(t_y).$$
	
Denote by $K_i$ the compact part of the level curve $\Sigma^+_i$ bounded by  $R_i$ and denote by  $V_i^3\subset U_{p_i}$ the compact set bounded by $R_i,K_i,d_i$. Define the function $\varphi_{_{V_i^3}}:V_i^3\to[i-2/3,i+1/3]$ by $\varphi_{_{V_i^3}}=\varphi_{p_i,i}.$
	
	Let $C_{i-1}=\partial U_{i-1}\setminus S_{i-1}$ and $V^4_{i}=\bigcup\limits_{x\in C_{i-1}}\left(\bigcup\limits_{t\in[0,1]}f^{-t}(x)\right)$. 
	
	Thus every point $y\in V_i^4$ is uniquely represented as $y=f^{-t_y}(x_y)$ where $x_y\in C_{i-1}$ and $t_y\in[0,1]$. Define the function $\varphi_{_{V_i^4}}:V_i^4\to[i-2/3,i+1/3]$ by $$\varphi_{_{V_i^4}}(y)=i+t_y-2/3.$$ 
	
	Let ${U}_{i}=U_{i-1}\cup\bigcup\limits_{i=1}^4 V_i^j$ and define the required function $\varphi_i$ by $$\varphi_{i}(x)=
	\begin{cases}
	{\varphi}_{i-1},  x\in U_{i-1}, \cr
	{\varphi}_{V_i^j},  x\in V_i^j,j\in\{1,2,3,4\}, 
	\end{cases}$$
	and that concludes the construction in this case.
	
{\it c) The point $p_{i}$ is a source.} From Proposition \ref{action_on_W} it follows that there exists a local energy function ${\varphi}_{p_{i},i}$ in some neighborhood of $p_{i}$. Let $\tilde{U}_{{i}}={\varphi}^{-1}_{p_{i},i}([i-1/3,i])$ and $\Sigma_{i}=\varphi^{-1}_{p_{i}, i}(i-1/3)$. Without loss of generality assume $\Sigma_i\cap U_{i-1} =\emptyset$ (to satisfy this condition one changes values of the function $\varphi_{p_{i},i}$). 
	
\begin{figure}[ht]
		\centering
	\includegraphics[width=0.7\linewidth]{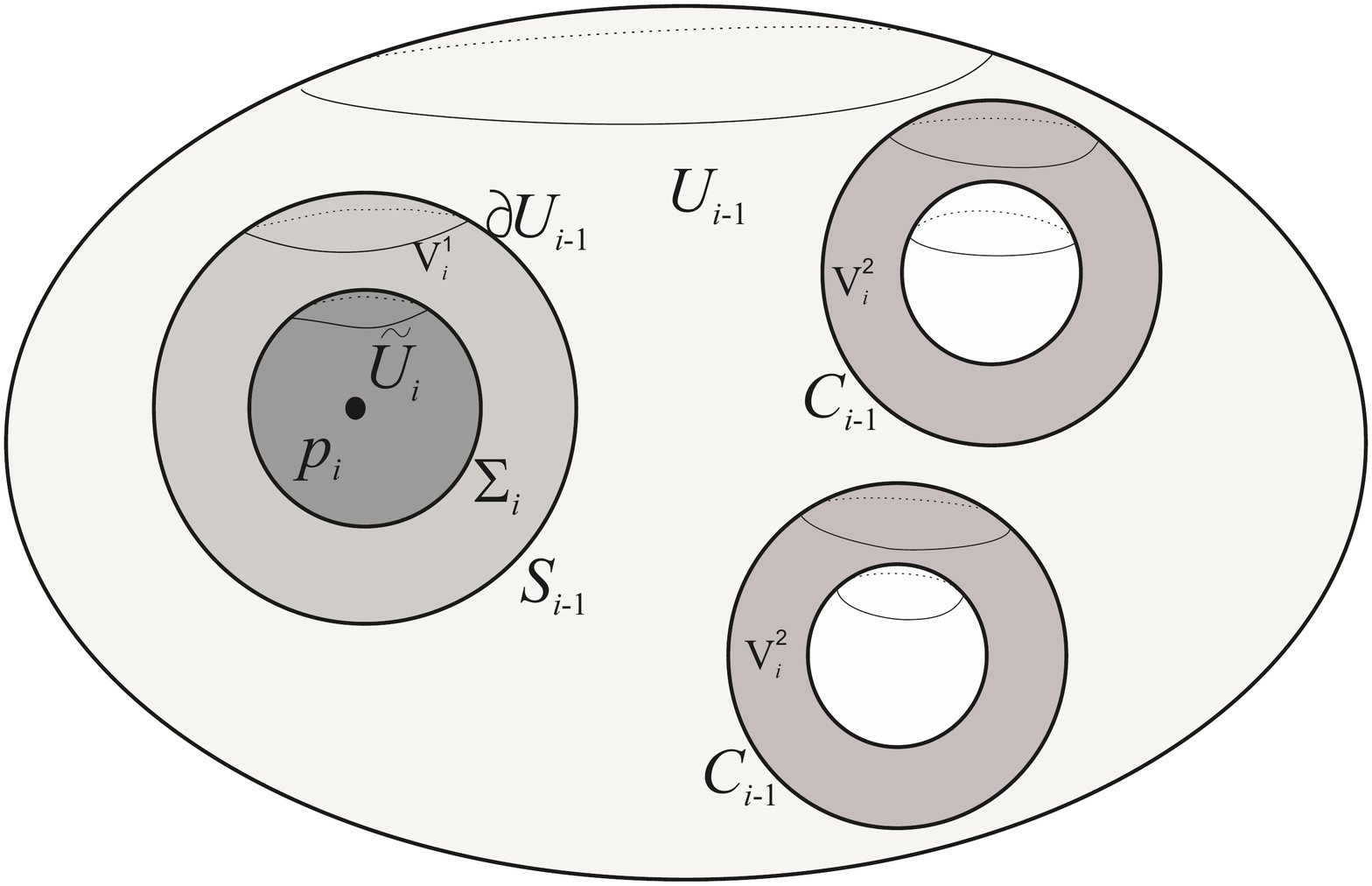}
		\caption{Induction step if $p_i$ is a source}
		\label{fig4}
	\end{figure}

Let $S_{i-1}$ be one of the connected components of $\partial U_{i-1}=\varphi^{-1}_{i-1}(i-2/3)$ and such that $\mathcal O_x\cap S_{i-1}\neq\emptyset, x\in \Sigma_{i}$ (see Figure \ref{fig4}). Then $f^{t_1(x)}(x)=\mathcal O_x\cap S_{i-1}, x\in \Sigma_{i}$ defines a positive continuous function $t_1:\Sigma_{i}\to\mathbb R$. Let $$V^1_i=\bigcup\limits_{x\in \Sigma_{i}}\left(\bigcup\limits_{t\in[0,t_1(x)]}f^{t}(x)\right).$$ Every point $y\in V^1_i$ is uniquely represented as $y=f^{t_y}(x_y)$ for some point $x_y\in\Sigma_{i}$ and some time moment $0\leq t_y\leq t_1(x_y)$. Define the function $\varphi_{_{V_i^1}}:V_i^1\to[i-2/3,i-1/3]$ by $$\varphi_{_{V_i^1}}(y)=i-1/3\left(1+t_y/t_1(x_y)\right).$$
	
Let $C_{i-1}=\partial U_{i-1}\setminus S_{i-1}$ and $V^2_{i}=\bigcup\limits_{x\in C_{i-1}}\left(\bigcup\limits_{t\in[0,1]}f^{-t}(x)\right)$. 
	
	Thus every point $y\in V_i^2$ is uniquely represented as $y=f^{-t_y}(x_y)$ where $x_y\in C_{i-1}$ and $t_y\in[0,1]$. Define the function $\varphi_{_{V_i^2}}:V_i^2\to[i-2/3,i+1/3]$ by $$\varphi_{_{V_i^2}}(y)=i+t_y-2/3.$$ Let ${U}_{i}=U_{i-1}\cup\tilde{U}_i\cup\bigcup\limits_{i=1}^2 V_i^j$ and define the required function $\varphi_i$ by
	$$\varphi_{i}= 
	\begin{cases}
	{\varphi}_{i-1},  x\in U_{i-1},\cr  
	{\varphi}_{p_{i},i},  x\in \tilde{U}_{i},\cr  
	{\varphi}_{_{V_i^j}},  x\in V_i^j,\,j\in\{1,2\}
	\end{cases}$$ and that concludes the construction in this case.
	
Thus by induction we have constructed the function $\varphi_{k-1}:U_{k-1}\to[1,k-2/3]$ with the desired properties. We now turn to construction of the required energy function $\varphi_{k}:U_{k}\to[1,k]$.
	
From Proposition \ref{action_on_W} it follows that there exists a local energy function ${\varphi}_{p_{k},k}$ in some neighborhood of the source $p_{k}$.  Let $\tilde{U}_{{k}}={\varphi}^{-1}_{p_{k},k}([k-1/3,k])$ and $\Sigma_{k}=\varphi^{-1}_{p_{k}, k}(k-1/3)$. Without loss of generality assume $\tilde U_k \cap U_{k-1}=\emptyset$ (to satisfy this condition one changes values of the function $\varphi_{p_{k},k}$). Then $f^{t_1(x)}(x)=\mathcal O_x\cap\partial U_{k-1}, x\in \Sigma_{k}$ defines a positive continuous function $t_1:\Sigma_{k}\to\mathbb R$. Let $$V_k=\bigcup\limits_{x\in \Sigma_{k}}\left(\bigcup\limits_{t\in[0,t_1(x)]}f^{t}(x)\right).$$ Every point $y\in V_k$ is uniquely represented as $y=f^{t_y}(x_y)$ for the point $x_y\in\Sigma_{k}$ and the time $0\leq t_y\leq t_1(x_y)$. Define the function $\varphi_{_{V_k}}:V_k\to[k-2/3,k-1/3]$ by $$\varphi_{_{V_k}}(y)=k-1/3\left(1+t_y/t_1(x_y)\right).$$
	
	Let $U_{k}= U_{k-1}\cup V_k$ and define the required function by
	$$\varphi_{k}= 
	\begin{cases}
	\varphi_{k-1}(x),  x\in U_{k-1},\cr  
	\varphi_{_{V_k}}(x),  x\in V_k,
	\end{cases}$$ and this concludes the construction.
	
		\begin{figure}[ht]
		\centering
		\includegraphics[width=0.7\linewidth]{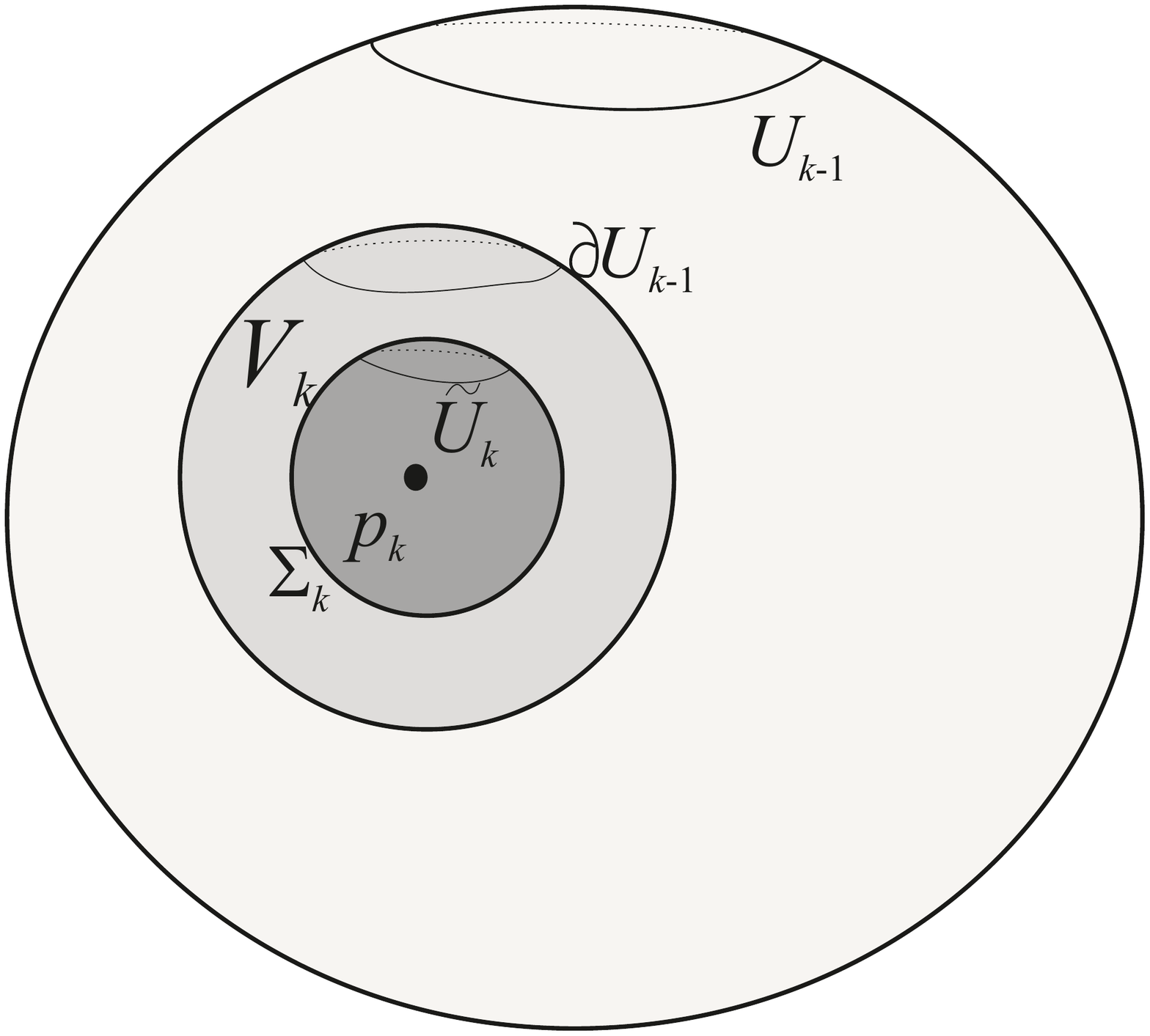}
		\caption{Construction of the required energy function $\varphi_{k}:U_{k}\to[1,k]$ }
		\label{fig5}
	\end{figure}

Funding: The study was implemented in the framework of the Basic Research Program at the National Research University Higher School of Economics (HSE University) in 2019.

\bibliography{literatshort}

\begin{thebibliography}{1}
\expandafter\ifx\csname url\endcsname\relax
  \def\url#1{\texttt{#1}}\fi
\expandafter\ifx\csname urlprefix\endcsname\relax\def\urlprefix{URL }\fi
\expandafter\ifx\csname href\endcsname\relax
  \def\href#1#2{#2} \def\path#1{#1}\fi

\bibitem{Con1978}
C.~C. Conley, Isolated invariant sets and the {Morse} index, Vol.~38, American
  Mathematical Soc., 1978.

\bibitem{wilson1969smoothing}
F.~Wilson, Smoothing derivatives of functions and applications, Trans. Amer.
  Math. Soc. 139 (1969) 413--428.

\bibitem{S12}
S.~Smale, On gradient dynamical systems, Ann. of Math. (1961) 199--206.

\bibitem{Me1968}
K.~R. Meyer, Energy functions for {Morse-Smale} systems, Amer. J. Math. (1968)
  1031--1040.

\bibitem{GrGuMePo2017}
V.~Grines, E.~Gurevich, V.~Medvedev, O.~Pochinka, An analog of {Smale's}
  theorem for homeomorphisms with regular dynamics, Math. Notes 102~(3-4)
  (2017) 569--574.

\bibitem{morse1959topologically}
M.~Morse, Topologically non-degenerate functions on a compact n-manifold m, J.
  Anal. Math. 7~(1) (1959) 189--208.

\bibitem{hir94}
M.~W. Hirsch, Differential Topology, Springer-Verlag, 1994.

\end{thebibliography}

\end{document}